\documentstyle{article}
\begin{document}
\begin{center}
{\bf SOLUTIONS OF FRACTIONAL REACTION-DIFFUSION EQUATIONS
IN TERMS OF THE H-FUNCTION}\par
\bigskip
H.J. HAUBOLD\\
Office for Outer Space Affairs, United Nations, Vienna International Centre\\
P.O. Box 500, A-1400, Vienna, Austria and\\
Centre for Mathematical Sciences, Pala Campus\\ 
Arunapuram P.O., Pala-686 574, Kerala, India\\[0.5cm]
A.M .MATHAI\\
Department of Mathematics and Statistics, McGill University\\
Montreal, Canada H3A 2K6 and \\
Centre for Mathematical Sciences, Pala Campus\\ 
Arunapuram P.O., Pala-686 574, Kerala, India\\[0.5cm]
R.K. SAXENA\\
Department of Mathematics and Statistics, Jai Narain Vyas University\\
Jodhpur-342004, India
\end{center}
\bigskip
\noindent
{\bf Abstract.} This paper deals with the investigation of the solution of an unified fractional reaction-diffusion equation associated with the Caputo derivative as the time-derivative and Riesz-Feller fractional derivative as the space-derivative. The solution is derived by the application of the Laplace and Fourier transforms in closed form in terms of the H-function. The results derived are of general nature and include the results investigated earlier by many authors, notably by Mainardi et al. (2001, 2005) for the fundamental solution of the space-time fractional diffusion equation, and Saxena et al. (2006a, b) for fractional  reaction- diffusion equations. The advantage of using Riesz-Feller derivative lies in the fact that the solution of the fractional reaction-diffusion equation containing this derivative includes the fundamental solution for space-time fractional diffusion, which itself is a generalization of  neutral fractional diffusion, space-fractional diffusion, and time-fractional diffusion. These specialized types of diffusion can be interpreted as spatial probability density functions evolving in time and are expressible in terms of the H-functions in compact form.
\section{Introduction}  
The review of the theory and applications of reaction-diffusion systems is contained in many books and articles. In recent work authors have demonstrated the depth of mathematics and related physical issues of reaction-diffusion equations such as nonlinear phenomena, stationary and spatio-temporal dissipative pattern formation, oscillations, waves etc. (Frank, 2005; Grafiychuk, Datsko, and Meleshko, 2006, 20076). In recent time, interest in fractional reaction-diffusion equations has increased because the equation exhibits self-organization phenomena and introduces a new parameter, the fractional index, into the equation. Additionally, the analysis of fractional reaction-diffusion equations is of great interest from the analytical and numerical point of view.

The objective of this paper is to derive the solution of an unified model of reaction-diffusion system (14), associated with the Caputo derivative and the Riesz-Feller derivative. This new model provides the extension of the models discussed earlier by Mainardi, Luchko, and Pagnini (2001), Mainardi, Pagnini, and Saxena (2005), and Saxena, Mathai, and Haubold (2006a). The present study is in continuation of our earlier work, Haubold and Mathai (1995, 2000) and  Saxena, Mathai, and Haubold (2006a, 2006b).

\section{Results Required in the Sequel}
In view of the results
\begin{equation}
J_{-1/2}(x)=\sqrt{\frac{2}{\pi x}}cos x.
\end{equation}
and (Mathai and Saxena, 1978, p. 49), the cosine transform of the H-function is given by 
\begin{eqnarray}
&&\int^\infty_0 t^{\rho-1} cos(kt)H^{m,n}_{p,q}\left[at^\mu\left|^{(a_p,A_p)}_{(b_q,B_q)}\right.\right] dt\\
&=&\frac{\pi}{k^\rho}\;H^{n+1,m}_{q+1,p+2}\left[\frac{k^\mu}{a}\left|^{(1-b_q, B_q), (\frac{1+\rho}{2},\frac{\mu}{2})}
_{(\rho, \mu), (1-a_p, a_p), (\frac{1+\rho}{2}, \frac{\mu}{2})}\right]\right.,
\end{eqnarray}
\noindent
where $Re[\rho+\mu^{min}_{1\leq j\leq m}(\frac{b_j}{B_j})]>0, Re[\rho+\mu^{max}_{1\leq j\leq n}\left(\frac{a_j-1}{A_j}\right)]<0, |arg \alpha|<\frac{1}{2}\pi \Omega, \Omega>0;\\
k> 0\; \mbox{and}\; \Omega = \sum^m_{j=1} B_j - \sum ^q_{j=m+1} B_j+\sum^n_{j=1} A_j-\sum^p_{j=n+1} A_j.$\\
      The Riemann-Liouville fractional integral of order $\nu$  is defined by 
 (Miller and Ross, 1993, p. 45; Kilbas et al., 2006)
\begin{equation}
_0D_t^{-\nu}N(x,t)=\frac{1}{\Gamma(\nu)}\int_0^t(t-u)^{\nu-1}N(x,u)du,
\end{equation}
where $Re(\nu)>0.$

The following fractional derivative of order $\alpha>0$   is introduced by Caputo (1969; see also Kilbas et al., 2006) in the form
\begin{eqnarray}
_0D_t^\alpha f(x,t)&=&\frac{1}{\Gamma(m-\alpha)}\int_0^t\frac{f^{(m)}(x,\tau)d\tau}{(t-\tau)^{\alpha+1-m}}, m-1<\alpha\leq m, Re(\alpha)>0, m\in N.\nonumber\\
&=& \frac{\partial^m f(x,t)}{\partial t^m},\; \mbox{if}\;\alpha=m.
\end{eqnarray}
where $\frac{\partial^m}{\partial^m} f(x,t)$  is the $m^{th}$ partial  derivative of f(x,t) with respect to t.

     The Laplace transform of the Caputo derivative is given by Caputo (1969; see also Kilbas et al., 2006) in the form
\begin{equation}
L\left\{_0D_t^\alpha f(x,t); s\right\}= s^\alpha F(x,s) - \sum^{m-1}_{r=0}s^{\alpha-r-1}f^{(r)}(x,0+),\; (m-1<\alpha\leq m).
\end{equation}  

     Following Feller (1952, 1971), it is conventional  to define  the Riesz-Feller space-fractional derivative of order $\alpha$ and skewness $\theta$   in terms of its Fourier transform as 
\begin{equation}
F\left\{_xD^\alpha_\theta f(x); k\right\} = -\Psi^\theta_\alpha (k) f^* (k),
\end{equation}                                                  
where
\begin{equation}
\Psi^\theta_\alpha(k)=|k|^\alpha exp[i(signk)\frac{\theta\pi}{2}], \;0<\alpha\leq 2, |\theta|\leq min \left\{\alpha, 2- \alpha\right\}.
\end{equation}
When $\theta=0$, then (8) reduces to 
\begin{equation}
F\left\{_xD_0^\alpha f(x); k\right\}= -|k|^\alpha,
\end{equation}              				         
which is the Fourier transform of  the Weyl fractional operator, defined by
\begin{equation}
_{-\infty}D_x^\mu f(t)=\frac{1}{\Gamma(n-\mu)}\frac{d^n}{dt^n}\int_{-\infty}^t\frac{f(u)du}{(t-u)^{\mu-n+1}}.
\end{equation}
    This shows that the Riesz-Feller operator may be regarded as a generalization of the Weyl operator.                                              

Further, when $\theta=0$, we have a symmetric operator with respect to $x$ that can be interpreted as
\begin{equation}
_xD_0^\alpha=-\left(-\frac{d^2}{dx^2}\right)^{\alpha/2}
\end{equation}
This can be formally deduced by writing $-(k)^
\alpha=-(k^2)^{\alpha/2}.$ For $0<\alpha<2$   and $|\theta|\leq min\left\{\alpha, 2-\alpha\right\}$, 
the Riesz-Feller derivative can be shown to possess the following integral representation in the $x$ domain:
\begin{eqnarray}
_xD_\theta^\alpha f(x)&=&\frac{\Gamma(1+\alpha)}{\pi}\left\{sin[(\alpha+\theta)\pi/2]\int_0^\infty\frac{f(x+\xi)-f(x)}{\xi^{1+\alpha}}d\xi \right.\nonumber\\
&+&\left.sin[(\alpha-\theta)\pi/2]\int_0^\infty\frac{f(x-\xi)-f(x)}{\xi^{1+\alpha}}d\xi\right\}.
\end{eqnarray}

                                                                                                                          Finally, we need the following property of the H-function (Mathai and Saxena, 1978)
\begin{equation}
H^{m,n}_{p,q}\left[x^\delta\left|^{(a_p, a_p)}_{(b_q, B_q)}\right.\right]= \frac{1}{\delta}H^{m,n}_{p,q}\left[x\left|^{(a_p, A_p/\delta)}_{(b_q, B_q/\delta}\right]\right.\;\;,(\delta>0).
\end{equation}

\section{Unified Fractional Reaction-Diffusion Equation}
    In this section, we will investigate the solution of the reaction-diffusion equation (14) under the initial conditions (15). The result is given in the form of the following \\ 
{\bf Theorem}.  Consider the unified fractional reaction-diffusion model 
\begin{equation}
_0D_t^\beta N(x,t)=\eta_xD_\theta^\alpha N(x,t)+\Phi(x,t),
\end{equation}                                                             
where $\eta, t>0, x\in r; \alpha, \theta, \beta$  are real parameters with the constraints\\
$0<\alpha\leq 2, |\theta|\leq min(\alpha, 2- \alpha), 0 <\beta \leq 2,$
and the initial conditions
\begin{equation}
N(x,0)=f(x), N_t(x,0)=g(x)\;);\mbox{for} \;\;x \in R, ^{lim}_{|x|\rightarrow \pm \infty} N(x,t) = 0, t>0.
\end{equation} 
Here $N_t(x,0)$   means the first partial derivative of $N(x, t)$ with respect to $t$  evaluated at $t=0, \eta$   is a diffusion constant  and $\Phi(x,t)$  is a nonlinear function belonging to the area of reaction-diffusion. Further $_xD_\theta ^\alpha$   is the Riesz-Feller space-fractional derivative of order $\alpha$  and asymmetry $\theta$. $_0D_t^\beta$  is the Caputo time-fractional  derivative  of order $\beta$. Then for the solution of (14), subject to the above constraints, there holds the formula 
\begin{eqnarray}
N(x,t)&=&\frac{1}{2\pi}\int^\infty_{-\infty}f^*(k)E_{\beta,1}(-\eta t^\beta \Psi^\theta_\alpha(k))exp(-ikx)dk\\
&+& \frac{1}{2\pi}\int^\infty_{-\infty} tg^*(k)E_{\beta,2}(-\eta k^\alpha t^\beta \Psi ^\theta_\alpha(k)) exp(-ikx)dk\nonumber\\
&+&\frac{1}{2\pi}\int_0^t \xi^{\beta-1}d\xi\int^\infty_{-\infty}\Phi^*(k,t-\xi)E_{\beta,\beta}(-\eta k^\alpha t^\beta\Psi^\theta_\alpha(k))exp(-ikx)dk.\nonumber
\end{eqnarray}
In equation (16) and the following, $E_{\alpha,\beta}(z)$ denotes the generalized Mittag-Leffler function (Saxena, Mathai, and Haubold, 2004; Berberan-Santos, 2005; Chamati and Tonchev, 2006).\\
{\bf Proof.} If we apply the Laplace transform with respect to the time variable $t$, Fourier transform with respect to space variable 
$x$, and use the initial conditions (15) and the formula (7), then the given equation transforms into the form
$$s^\beta N^{^*_\sim}(k,s)-s^{\beta-1}f^*(k)-s^{\beta-2}g^*(k)=-\eta\Psi^\theta_\alpha(k)N^{^*_\sim}(k,s)+\Phi^{^*_\sim}(k,s),$$
where according to the conventions followed , the symbol $\sim$ will stand for the Laplace transform with respect to time variable $t$ and *  represents the Fourier transform with respect to space variable $x$. 

Solving for $N^{^*_\sim}$, it yields
\begin{equation}
N^{^*_\sim}(k,s)=\frac{f^*(k)s^{\beta-1}}{s^\beta+\eta\Psi^\theta_\alpha(k)}+\frac{g^*(k)s^{\beta-2}}{s^\beta+\eta\Psi^\theta_\alpha(k)}+\frac{\Phi^{^*_\sim}(k)}{s^\beta+\eta\Psi^\theta_\alpha(k)}.
\end{equation}
On taking the inverse Laplace transform of (17) and applying the formula 
\begin{equation}
L^{-1}\left\{\frac{s^{\beta-1}}{a+s^\alpha}\right\}=t^{\alpha-\beta}E_{\alpha, \alpha-\beta+1}(-at^\alpha),
\end{equation}
where $Re(s) > 0, Re(\alpha)>0, Re(\alpha-\beta)>-1;$ it is seen that 
\begin{eqnarray}
N^*(k,t) &=& f^* (k) E_{\beta,1}(-\eta t^\beta\Psi_\alpha^\theta(k))+g^*(k)tE_{\beta, 2}(-\eta t^\beta \Psi^\theta_\alpha(k))\nonumber\\
&+&\int_0^t\Phi^*(k,t-\xi)\xi^{\beta-1}E_{\beta,\beta}(-\eta\Psi_\alpha^\theta(k)\xi^\beta)d\xi.
\end{eqnarray}
The required solution (16) is now obtained by taking the inverse Fourier transform of (19). This completes the proof of the theorem.

\section{Special Cases}
    When $g(x) = 0$, then by the application of the convolution theorem of the Fourier transform to the solution (16) of the theorem, it readily yields\\   
{\bf Corollary 1.} The solution of the fractional reaction-diffusion equation 
\begin{equation}
\frac{\partial^\beta}{\partial t^\beta}N(x,t) - \eta\frac{\partial^\alpha}{\partial x^\alpha} N(x,t) = \Phi(x,t), x\in R, t>0, \eta>0,
\end{equation}
with initial conditions  
\begin{equation}
N(x,0)=f(x), N_t(x,0)=0 \;\mbox{for}\; x\in R, 1<\beta\leq 2, ^{lim}_{x\rightarrow \pm \infty} N(x,t)=0,
\end{equation}
where $\eta$ is  a diffusion constant and $\Phi(x,t)$ is a nonlinear function belonging to the area of reaction-diffusion, is given by 
\begin{eqnarray}
    N(x,t)  &=& \int_0^x G_1(x-\tau, t) f(\tau) d\tau\nonumber\\
&+&\int_0^t(t-\xi)^{\beta-1}d\xi\;\;\int_0^x G_2(x-\tau, t-\xi)\Phi(\tau, \xi) d\tau,
\end{eqnarray}  

where
\begin{eqnarray}
\rho&=&\frac{\alpha-\theta}{2\alpha}\nonumber\\
G_1(x,t)&=&\frac{1}{2\pi}\int_{-\infty}^{\infty} exp(-ikx)E_{\beta, 1}(-\eta|t^\beta|\Psi^\theta_\alpha(k))dk\\
&=&\frac{1}{\alpha|x|}H^{2,1}_{3,3}\left[\frac{|x|}{\eta^{1/\alpha} t^{\beta/\alpha}}\left|^{(1,1/\alpha),(\beta, \beta/\alpha), (1,\rho)}_{(1,1/\alpha), (1,1), (1,\rho)}\right.\right], (\alpha>0)\nonumber 
\end{eqnarray}
and  
\begin{eqnarray}
G_2(x,t)&=& \frac{1}{2\pi}\int^\infty_{-\infty}exp(-ikx)E_{\beta, \beta}(-\eta t^\beta \Psi^\theta_\alpha(k))dk\nonumber\\
&=&\frac{1}{\alpha|x|}H^{2,1}_{3,3}\left[\frac{|x|}{\eta^{1/\alpha}t^{\beta/\alpha}}\left|^{(1,1/\alpha), (\beta, \beta/\alpha),(1,\rho)}_{(1,1/\alpha),(1,1), (1,\rho)}\right.\right]\;,(\alpha>0).
\end{eqnarray}
In deriving the above results, we have used the inverse Fourier transform formula 
\begin{equation}
F^{-1}[E_{\beta,\gamma}(-\eta t^\beta \Psi_\theta^\alpha(k));x]=\frac{1}{\alpha|x|}H^{2,1}_{3,3}[\frac{|x|}{\eta^{1\alpha}t^{\beta/\alpha}}|^{(1,1/\alpha),(\gamma,\beta/\alpha),(1,\rho)}_{(1,1/\alpha),(1,1),(1,\rho)}],
\end{equation}
where $Re(\beta)>0, Re(\gamma)>0,$ which can be established by following a procedure similar to that employed by Mainardi, Luchko, and Pagnini (2001).
     Next , if we set  $f(x)=\delta(x), \Phi=0, g(x)=0,$ where $\delta(x)$ is the Dirac delta-function, then we arrive at the following interesting result given by Mainardi, Pagnini, and Saxena (2005).\\
{\bf Corollary 2.} Consider the following space-time fractional diffusion model 
\begin{equation}
\frac{\partial^\beta N(x,t)}{\partial t^\beta}=\eta\;_xD_\theta^\alpha N(x,t), \eta >0, x\in R, \;\;0<\beta\leq 2,
\end{equation}
with the initial conditions  $N(x,t = 0) = \delta(x), N_t(x,0)=0, ^{lim}_{x\rightarrow \pm\infty}N(x,t)=0$ where $\eta$ is a diffusion constant  and $\delta(x)$  is the Dirac delta-function. Then for the fundamental  solution of (26)  with initial conditions, there holds  the formula 
\begin{equation}
N(x,t)=\frac{1}{\alpha|x|}H^{2,1}_{3.3}[\frac{|x|}{(\eta t^\beta)^{1/\alpha}}|^{(1,1/\alpha), (1,\beta/\alpha), (1,\rho)}_{(1,1/\alpha), (1,1),(1,\rho)}],
\end{equation}
where $\rho=\frac{\alpha - \theta}{2\alpha}.$\\ 
Some interesting special cases of (26) are enumerated below.

    (i)  We note that for $\alpha=\beta$, Mainardi, Pagnini, and Saxena (2005) have shown that the corresponding solution of (26), denoted by $N^\theta_\alpha$, which we call as the neutral fractional diffusion, can be expressed in terms of elementary function  and  can be defined for $x>0$ as\\
Neutral fractional diffusion: $0<\alpha=\beta<2; \theta\leq min\left\{\alpha, 2-\alpha\right\},$
\begin{equation}
N^\theta_\alpha(x)=\frac{1}{\pi}\frac{x^{\alpha-1} sin[(\pi/2)(\alpha-\theta)]}{1+2x^\alpha cos[(\pi/2)(\alpha-\theta)]+x^{2\alpha}}.
\end{equation}
The neutral fractional diffusion is not studied at length in the literature. 

      Next we derive some stable densities in terms of the H-functions as special cases of the solution of the equation (26)

       (ii) If we set $\beta=1,0<\alpha<2; \theta \leq min\left\{\alpha, 2 -\alpha \right\}$then (26) reduces to space fractional diffusion equation, which we denote by $L^\theta_\alpha (x)$ is the fundamental solution of the following  space-time fractional diffusion model:
\begin{equation}
\frac{\partial N(x,t)}{\partial t} = \eta\;_xD_\theta^\alpha N(x,t), \;\;\eta>0, \;x\in R,
\end{equation}
with the initial conditions  $N (x,t = 0) =\delta(x),\;^{lim}_{x\rightarrow\pm \infty} N(x,t)=0,$, where $\eta$ is a diffusion constant  and $\delta(x)$ is the Dirac-delta function. Hence for the solution of (29) there holds the formula 

\begin{equation}
L^\theta_\alpha (x)=\frac{1}{\alpha(\eta t)^{1/\alpha}}H^{1,1}_{2,2}\left[\frac{(\eta t)^{1/\alpha}}{|x|}\left|^{(1,1),(\rho,\rho)}_{(\frac{1}{\alpha},\frac{1}{\alpha}), (\rho, \rho)}\right.\right], \; 0<\alpha<1, |\theta|\leq \alpha,
\end{equation}
where $\rho=\frac{\alpha-\theta}{2\alpha}$. The density represented by the above expression is known as $\alpha$-stable L\'{e}vy density. Another form of this density is given by 
\begin{equation}
L^\theta_\alpha (x)=\frac{1}{\alpha(\eta t)^{1/\alpha}}H^{1,1}_{2,2}\left[\frac{|x|}{(\eta t)^{1/\alpha}}\left|^{(1-\frac{1}{\alpha},\frac{1}{\alpha}),(1-\rho, \rho)}_{(0,1), (1-\rho, \rho)}\right], \;1<\alpha<2, |\theta|\leq 2-\alpha,\right.
\end{equation}

(iii) Next, if we take $\alpha = 2,0 <\beta < 2,\theta =0,$  then we obtain the time fractional diffusion,  which is governed by the following  time fractional diffusion model:
\begin{equation}
\frac{\partial^\beta N(x,t)}{\partial t^\beta} = \eta\frac{\partial ^2}{\partial x^2}N(x,t), \eta>0, x\in R, 0<\beta \leq 2,
\end{equation}
with the initial conditions  $N (x,t = 0) = \delta(x), N_t(x,0)=0, ^{lim}_{x\rightarrow \pm \infty}N(x,t)=0$ where $\eta$ is a diffusion constant  and $\delta(x)$ is the Dirac delta-function, whose fundamental  solution is given by the equation 
\begin{equation}
N(x,t)= \frac{1}{2|x|}H^{1,0}_{1,1}\left[\frac{|x|}{(\eta t^\beta)^{1/2}}\left|^{(1,\beta/2)}_{(1,1)}\right.\right].
\end{equation}

(iv) Further, if we set $\alpha=2, \beta=1$ and $\theta \rightarrow 0$  then for the fundamental solution of the standard diffusion equation
\begin{equation}
\frac{\partial}{\partial t}N(x,t)=\eta\frac{\partial ^2}{\partial x^2} N(x,t),
\end{equation}
 with initial condition 
\begin{equation}
        N(x,t=0) = \delta(x),\;^{lim}_{x\rightarrow\pm\infty} N(x,t) =0,
\end{equation}
there holds  the formula 
\begin{equation}
N(x,t)=\frac{1}{2|x|}H^{1,0}_{1,1}\left[\frac{|x|}{\eta^{1/2}t^{1/2}}\left|^{(1,1/2)}_{(1,1)}\right]=(4\pi\eta t)^{-1/2}exp[-\frac{|x|^2}{4\eta t}],\right.
\end{equation}               
which is the classical Gaussian density. For further details of these special cases based on the Green function, one can refer to the paper by Mainardi, Luchko, and Pagnini (2001) and Mainardi, Pagnini, and Saxena (2005).\\
{\bf Remark.} Fractional order moments and the asymptotic expansion of the solution (27) are discussed by Mainardi, Luchko, and Pagnini (2001).
        
      Finally, for $\beta=1/2$ in (14), we arrive at\\ 
{\bf Corollary 3.} Consider the following fractional reaction-diffusion model 
\begin{equation}
D_t^{1/2}N(x,t)=\eta_x D_\theta^\alpha N(x,t)+\Phi(x,t),
\end{equation}
where $\eta, t >0, x\in R; \alpha, \theta$  are real parameters with the constraints\\
$0<\alpha\leq 2, |\theta|\leq min(\alpha, 2-\alpha),$    
and the initial conditions 
\begin{equation}
N(x,0)= f(x),\; \mbox{for}\;x\in R, \;^{lim}_{x\rightarrow\pm \infty} N(x,t)=0.
\end{equation}
Here $\eta$   is a diffusion constant and $\Phi(x,t)$  is a nonlinear function belonging to the area of reaction-diffusion. Further $_xD_\theta^\alpha$  is the Riesz-Feller space fractional derivative of order $\alpha$  and asymmetry $\theta$  and $D_t^{1/2}$ is the Caputo time-fractional derivative of order $1/2$. Then for the solution of (37), subject to the above constraints, there holds the formula 
\begin{eqnarray}
N(x,t)&=&\frac{1}{2\pi}\int^\infty_{-\infty} f^*(k)E_{1/2,1}(-\eta t^\beta \Psi^\theta_\alpha(k))exp(-ikx)dk\\
&+& \frac{1}{2\pi}\int_0^t\xi^{-1/2}d\xi\int^\infty_{-\infty} \Phi^*(kct-\xi)E_{\frac{1}{2},\frac{1}{2}}(-\eta k^\alpha t^{1.2}\Psi^\theta_\alpha(k))exp(-ikx)dk.\nonumber
\end{eqnarray}
If we set $\theta=0$  in (39), then it reduces to the result recently obtained by the authors (2006a) for the fractional reaction-diffusion equation.\par
\bigskip
\noindent
\section {References}\par
\medskip
\noindent
Berberan-Santos, M.N. (2005). Properties of the Mittag-Leffler relaxation function, {\it Journal of Mathematical Chemistry,} {\bf 38}, 629-635.\par
\medskip
\noindent
Caputo, M. (1969). {\it Elasticita  e Dissipazione}, Zanichelli, Bologna.\par
\medskip
\noindent
Chamati, H. and Tonchev, N.S. (2006). Generalized Mittag-Leffler functions in the theory of finite-size scaling for systems with strong anisotropy and/or long-range interaction, {\it Journal of Physics A: Mathematical and General,} {\bf 39}, 469-478.\par
\medskip
\noindent
Feller, W. (1952). On a generalization of Marcel Riesz' potentials and the semi-groups generated by them, {\it Meddeladen Lund Universitets Matematiska Seminarium }(Comm. S\'{e}m. Math\'{e}m. Universit\'{e} de Lund ), Tome suppl. d\'{e}di\'{e} a M. Riesz, Lund, 73-81.\par
\medskip
\noindent
Feller, W. (1966). {\it An Introduction to Probability Theory and its Applications}, Vol. II, John Wiley and Sons, New York.\par
\medskip
\noindent
Frank, T.D. (2005). {\it Nonlinear Fokker-Planck Equations: Fundamentals and Applications}, Springer, Berlin Heidelberg New York.\par
\medskip
\noindent
Grafiychuk, V., Datsko, B., and Meleshko, V. (2006). Mathematical modeling of pattern formation in sub- and superdiffusive reaction-diffusion systems, arXiv:nlin.AO/06110005 v3.\par
\medskip
\noindent
Grafiychuk, V., Datsko, B., and Meleshko, V. (2007). Nonlinear oscillations and stability domains in fractional reaction-diffusion systems, arXiv:nlin.PS/0702013 v1.\par
\medskip
\noindent
Haubold, H.J. and Mathai, A.M. (2000). The fractional kinetic equation and thermonuclear functions, {\it Astrophysics and Space Science,} {\bf 273}, 53-63.\par
\medskip
\noindent
Haubold, H.J. and Mathai, A.M. (1995). A heuristic remark on the periodic variation in the number of solar neutrinos detected on Earth, {\it Astrophysics and Space Science,} {\bf 228}, 113-124.\par
\medskip
\noindent
Kilbas, A.A., Srivastava, H.M., and Trujillo, J.J. (2006). {\it Theory and Applications of Fractional Differential Equations}, Elsevier, Amsterdam.\par
\medskip
\noindent 
Mainardi, F., Luchko, Y., and Pagnini, G. (2001). The fundamental solution of the space-time fractional diffusion equation, {\it Fractional Calculus and Applied Analysis.} {\bf 4}, 153-192.\par
\medskip
\noindent 
Mainardi, F., Pagnini, G., and Saxena, R.K. (2005). Fox H-functions in fractional diffusion, {\it Journal of Computational and Applied Mathematics} {\bf 178}, 321-331.\par
\medskip
\noindent 
Mathai, A.M. and Saxena, R.K. (1978). {\it The H-function with Applications in Statistics and Other Disciplines}, John Wiley and Sons, New York, London, and Sydney.\par
\medskip
\noindent 
Miller, K.S. and Ross, B. (1993). {\it An Introduction to the Fractional Calculus and Fractional Differential Equations}, John Wiley and Sons, New York.\par
\medskip
\noindent
Saxena, R.K., Mathai, A.M., and Haubold, H.J. (2004). On fractional kinetic equations, {\it Astrophysics and Space Science,} {\bf 282}, 281-287.\par
\medskip
\noindent 
Saxena, R.K., Mathai, A.M., and Haubold, H.J. (2006a). Fractional reaction-diffusion equations, {\it Astrophysics and Space Science,} {\bf 305}, 289-296.\par
\medskip
\noindent 
Saxena, R.K., Mathai, A.M., and Haubold, H.J. (2006b). Reaction-diffusion systems and nonlinear waves, Astrophysics and Space Science, {\bf 305}, 297-303.\par
\medskip
\noindent 
Yu, R. and Zhang, H. (2006). New function of Mittag-Leffler type and its application in the fractional diffusion-wave equation, Chaos, Solitons and Fractals {\bf 30}, 946-955.\par
\end{document}